\documentclass[12pt,oneside,final]{amsart}
\usepackage[cp1250]{inputenc}
\usepackage[T1]{fontenc}
\usepackage{amssymb}
\usepackage{color}
\renewcommand{\thefootnote}{}
\newtheorem{twr}{Twierdzenie}%[section]
\newtheorem{thr}[twr]{Theorem}

\newtheorem{lm}[twr]{Lemma}

\newtheorem{crl}[twr]{Corollary}
\renewcommand{\thefootnote}{}

\newcommand{\cn}{\mathbb{C}^n}

\newcommand{\dc}{i{\partial}\bar\partial}

\title[$J$-plurisubharmonic functions]{On regularization of $J$-plurisubharmonic functions  }
\author[S. Pli\'{s}]{Szymon Pli\'{s}}
\address{  Institute of Mathematics, Cracow University of Technology, Warszawska 24, 31-155
    Krak\'{o}w, Poland
}
\email{splis@pk.edu.pl}

\begin{document}\thispagestyle{empty} \footnotetext{The  author was partially supported by the NCN grant 2011/01/D/ST1/04192.
}\renewcommand{\thefootnote}{\arabic{footnote}}

\begin{abstract}
We show that on almost complex surfaces plurisubharmonic functions can be locally approximated by smooth plurisubharmonic functions. The main tool is the Poletsky type theorem due to U. Kuzman.
\end{abstract}

\maketitle
\section{Introduction}
 Let $(M,J)$ be an almost complex manifold. In his paper \cite{h} Haggui defines plurisubharmonic functions on $M$ as upper semicontinuous functions which are subharmonic on every $J$-holomorphic dysk. Recently Harvey and Lawson proved that a locally integrable function $u$ is plurisubharmonic iff a current $i\partial\bar{\partial}u$ is positive (see \cite{h-l1}).

 It is a very natural open question in this theory  whether any plurisubharmonic function is (locally) a limit of a decreasing sequence of smooth plurisubharmonic functions.  The Richberg type theorem was proved in \cite{p}. This gives a positive answer in a case of continuous functions. In this note we prove it for all plurisubharmonic functions in the (complex) dimension\footnote{In this note by the dimension of an almost complex manifold we mean the complex dimension which is a half of the real dimension.} 2.

\begin{thr}\label{aproksymacjaciaglymi} Let $\dim M=2$ and $P\in M$. Then there is a domain $D$ which is a neighbourhood of $P$ such that 
for every  $u\in\mathcal{PSH}(D)$ there exists a decreasing sequence          $\psi_k\in\mathcal{C}^\infty\cap\mathcal{PSH}(D)$ such that $\psi_k\rightarrow u$.
\end{thr}

As an immediate consequence of  Theorem \ref{aproksymacjaciaglymi} and  proposition 5.2 from  \cite{p} we obtain the following
\begin{crl}
Let $\dim M=2$ and $u,v$ in $W^{1,2}_{loc}\cap\mathcal{PSH}(M)$. Then a current $\dc u\wedge\dc v$ defined in \cite{p} is a (positive) measure.
\end{crl}

In particular   the Monge-Amp\`ere operator $(\dc u)^2$ is well defined for any bounded psh function $u$ on an almost complex surface (compare proposition 4.2  in \cite{p}). On domains in $\mathbb{C}^2$ it was proved by B\l ocki (see \cite{b}) that a set $W^{1,2}_{loc}\cap\mathcal{PSH}$ is a natural domain for the Monge-Amp\`ere operator.

The main step in the prove of Theorem \ref{aproksymacjaciaglymi} is the continuity of  largest plurisubharmonic minorants of  certain continuous functions. Harvey and Lawson, after viewing a preliminary version of this paper, inform me that using viscositi methods (as in section 8.6 in \cite{h-l2}) it is possible to prove it in any dimension. Let us stress that their methods are very different from methods used in this note.

\section{Proof}
\subsection{$J$-holomorphic discs}

A good reference for the (local) theory of $J$-holomorphic discs is \cite{i-r}.  In this subsection $J$ is $\mathcal{C}^{1}$ close to $J_{st}$ (in particular $(J+J_{st})$ is invertible) where $J_{st}$ is the standard (integrable) almost complex structure in $\cn$. Let $\mathbb{D}$ be a unit disc in $\mathbb{C}$. A function $u:\mathbb{D}\rightarrow (\cn,J)$ is $J$-holomorphic if and only if $$\frac{\partial u	}{\partial\bar z}+Q(u)\frac{\partial u	}{\partial z}=0 $$
when $$Q=(J-J_{st})(J+J_{st})^{-1}.$$
Let $0<\alpha<1$ and $T:\mathcal{C}^{0,\alpha}(\bar{\mathbb{D}},\cn)\rightarrow\mathcal{C}^{1,\alpha}(\bar{\mathbb{D}},\cn)$ be the Cauchy-Green operator given by
$$Tu=\frac{1}{\pi}\int_{\mathbb{D}}\frac{u(\zeta)}{\cdot-\zeta}d\zeta\;.$$
Set $$\Phi u=u+T(Q(u)\frac{\partial u	}{\partial z})$$
and $$\Psi u=\Phi u+(u-\Phi u)(0)\;.$$
By the definition $\Psi u(0)=u(0)$. Note that $u\in\mathcal{C}^{1,\alpha}(\bar{\mathbb{D}},\cn)$ is $J$-holomorphic in $\mathbb{D}$ iff $\Phi(u)$ is $J_{st}$-holomorphic.
Because $d\Psi$ is close to $Id$, the map $\Psi:\mathcal{C}^{1,\alpha}(\bar{\mathbb{D}},\cn)\rightarrow\mathcal{C}^{1,\alpha}(\bar{\mathbb{D}},\cn)$ is a local diffeomorphism and there is a constant $C_0$ such that $\|(d\Psi)^{-1}\|\leq C_0$ everywhere.

We will use the following 
\begin{lm}\label{przesuwanie dysku}
 Let $V\in\cn$. For any $u\in\mathcal{C}^{1,\alpha}(\bar{\mathbb{D}},\cn)$ there is  $v\in\mathcal{C}^{1,\alpha}(\bar{\mathbb{D}},\cn)$ such that $\Psi(v)=\Psi(u)+V$ and $\|u-v\|_{\mathcal{C}^{1,\alpha}}\leq C_0|V|$.
\end{lm}

\textit{Proof:} Set $U_t=\Psi(u)+tV$ and $$S=\\ \{t\in[0,1]: \exists w\in\mathcal{C}^{1,\alpha}(\bar{\mathbb{D}},\cn) \hbox{ s. t. }\Psi(w)=U_t,\;\|u-w\|_{\mathcal{C}^{1,\alpha}}\leq tC_0|V|\} .$$
$S$ is nonempty, by the inverse function theorem it is open, by the Arzel\`{a}-Ascoli theorem it is closed and hence $S=[0,1]$. \;$\Box$

\subsection{Disc envelope}

Let $p\in\Omega\subset M$ and let $\mathcal{O}_p(\bar{\mathbb{D}},\Omega)$ be a set of $J$-holomorphic discs $\lambda:\bar{\mathbb{D}}\rightarrow\Omega$ with $\lambda(0)=p$. For an upper semicontinuous function $f:\Omega\rightarrow\mathbb{R}$  we consider the following disc envelope:
$$P_\Omega f(p) =\inf_{\lambda\in\mathcal{O}_p(\bar{\mathbb{D}},\Omega)}\frac{1}{2\pi}\int_0^{2\pi}f\circ\lambda(e^{it})dt .$$

We need the following Lemma.
\begin{lm}\label{ciaglosc}
 Let $\Omega_1\Subset\Omega_2\subset\cn$ and let $J$ be an almost complex structure on $\Omega_2$ which is $\mathcal{C}^{1}$ close to $J_{st}$. Let $f\in\mathcal{C}(\Omega_2)$ be such that $$P_{\Omega_1}f=(P_{\Omega_2}f)|_{\Omega_1}.$$
Then $P_{\Omega_1}f\in\mathcal{C}(\Omega_1)$.
\end{lm}

\textit{Proof:} We can assume that $f$ is uniformly continuous on $\Omega_2$ with  a modulus of continuity $\omega$ and $J$ is $\mathcal{C}^{1}$ close to $J_{st}$ on $\cn$. Set any $0<\delta<C_0^{-1}{\rm dist}(\partial\Omega_1,\partial\Omega_2)$. Let $\varepsilon>0$, and $p,q\in\Omega_1$ with $|p-q|\leq\delta$. There is $\lambda\in\mathcal{O}_p(\bar{\mathbb{D}},\Omega_1)$ such that: $$P_{\Omega_1}f(p)\geq\frac{1}{2\pi}\int_0^{2\pi}f\circ\lambda(e^{it})dt-\varepsilon.$$
By Lemma \ref{przesuwanie dysku} there is a function $\mu\in\mathcal{C}^{1,\alpha}(\bar{\mathbb{D}},\cn)$ such that $\Phi(\mu)=\Phi(\lambda)+w-z$ and $\|\lambda-\mu\|_{L^\infty}\leq C_0|p-q|.$ Since functions $(z\mapsto \mu(rz))$ are in $\mathcal{O}_q(\bar{\mathbb{D}},\Omega_2)$ for $1>r>0$ we can estimate
$$P_{\Omega_1}f(q)=P_{\Omega_2}f(q)\leq\frac{1}{2\pi}\int_0^{2\pi}f\circ\mu(e^{it})dt$$ $$\leq\frac{1}{2\pi}\int_0^{2\pi}(f\circ\lambda(e^{it})+\omega(C_0\delta))dt\leq P_{\Omega_1}f(p)+\omega(C_0\delta)+\varepsilon.$$ Letting $\varepsilon$ to $0$ we can conclude that $P_{\Omega_1}f$ is uniformly continuous with  a modulus of continuity $\tilde\omega(x)=\omega(C_0x)$.\;$\Box$

\subsection{Kuzman-Poletsky theorem}
For a domain $\Omega\subset M=\cn$ and an upper semicontinuous function $f$, Poletsky (see \cite{po}) proved that $H_\Omega f$ is a plurisubharmonic function (moreover it is the largest plurisubharmonic minorant of $f$). The key tool in the proof of Theorem \ref{aproksymacjaciaglymi} is a result of Kuzman, who showed the same for any 2-dimensional almost complex manifold (see theorem 1 in \cite{k}). The only reason for the assumption about a dimension in our Theorem is  just this assumption in  Kuzman's theorem.

\textit{Proof of Theorem \ref{aproksymacjaciaglymi}:}
The theorem is local hence we can assume that $P\in\mathbb{C}^2$ and $J$ is $\mathcal{C}^{1}$ close to $J_{st}$. We can choose  a neighbourhood $D$ of $P$ such that there exists a positive continuous  strictly $J$-plurisubharmonic\footnote{A function $u$ is strictly plurisubharmonic on $D$ means as usually that for any $\varphi\in\mathcal{C}^2_0$ there is $\varepsilon>0$ such that $u+\varepsilon\varphi$ is plurisubharmonic. We write here $J$-plurisubharmonic instead of plurisubharmonic to stress that a function is plurisubharmonic with respect to the almost complex structure $J$ (note that on $D$ we have also the almost complex structure $J_{st}$).}  exhaustion function $\rho$  on $D$. \footnote{Such  domain $D$ is called a Stain domain, see \cite{d-s}. Here we can take $D=\{|z-P|<\varepsilon\}$ for $\varepsilon>0$ small enough.} Set $u\in\mathcal{PSH}(D)$. Let us take a decreasing sequence of continuous functions $\phi_k$ tending to $u$.  We can modify $\rho$ such that $\lim_{z\rightarrow\partial D}(\rho-\phi_1)=+\infty$ and put $\tilde{\phi}_k=\max\{\phi_k,\rho-k\}$. There are domains $D_k\Subset D$ such that  $\tilde{\phi}_k=\rho-k$ on some neighbourhood $U_k$ of $D\setminus D_k$. By Kuzman's result $\hat{\phi}_k=P_D\tilde{\phi}_k,P_{D_k}\tilde{\phi}_k$ are $J$-plurisubharmonic. Note that $\hat{\phi}_k=\rho-k$ on $U_k$ and $P_{D_k}\tilde{\phi}_k=\rho-k$ on $D_k\cap U_k$, hence by Lemma \ref{ciaglosc} $\hat{\phi}_k\in\mathcal{C}(D)$. Thus we get a decreasing sequence of continuous $J$-plurisubharmonic functions $\hat\phi_k$ tending to $u$.

By the Richberg theorem (see theorem 3.1 in \cite{p}) there are  functions $\psi_k\in\mathcal{C}^\infty\cap\mathcal{PSH}(D)$ such that $$\hat\phi_k+2^{-k-1}\rho\leq\psi_k\leq \hat\phi_k+2^{-k}\rho$$
and  we can see that a sequence $\psi_k$ decreases to u. $\Box$

\textbf{Acknowledgments.} The author would like to express his
gratitude to A. Sukhov for helpful discussions  
 on the subject of this paper.

\end{document}